\documentclass[12pt]{article}
\usepackage{amssymb}
\usepackage{amsmath}
\usepackage{amsfonts}
\usepackage{graphicx}
\usepackage{graphics,longtable}

\begin{document}

\begin{center}
\textbf{A new algorithm that generates the image of the attractor of a
generalized iterated function system}

\bigskip

\textit{Radu MICULESCU}, \textit{Alexandru MIHAIL} and \textit{%
Silviu-Aurelian URZICEANU}

\bigskip

\textbf{Abstract}

\bigskip
\end{center}

{\small We provide a new algorithm (called the grid algorithm) designed to
generate the image of the attractor of a generalized iterated function
system on a finite dimensional space and we compare it with the
deterministic algorithm regarding generalized iterated function systems
presented by P. Jaros, \L . Ma\'{s}lanka and F. Strobin in [Algorithms
generating images of attractors of generalized iterated function systems,
Numer. Algorithms, 73 (2016), 477-499].}

\bigskip



\bigskip

\textbf{I. Introduction}

\bigskip

As part of the effort to extend the classical theory of iterated function
systems due to J. Hutchinson (see [2]), R. Miculescu and A. Mihail
introduced the concept of generalized iterated function system (see [7] and
[9]) which was obtained by considering contractions from $X^{p}$ into $X$
rather than contractions from $X$ into itself (here $(X,d)$ is a metric
space and $p$ is a natural number). Sufficient conditions for the existence
and uniqueness of the attractor of a generalized iterated function system
(for short GIFS) $\mathcal{F}=((X,d),(f_{i})_{i\in \{1,2,...,L\}})$, an
upper bound for the Hausdorff--Pompeiu distance between the attractors of
two such GIFSs, an upper bound for the Hausdorff--Pompeiu distance between
the attractor of such a GIFS and an arbitrary compact set of $X$ have been
provided and the continuous dependence of the attractor on the functions $%
f_{i}$ was proved. In the last years this concept has been intensively
studied. Let us mention some lines of research regarding this subject:

In [15], F. Strobin and J. Swaczyna extended the concept of GIFS by using
weaker types of generalized contractions which are similar to those
introduced by F. Browder (see [1]) or J. Matkowski (see [5]). In [14],
Strobin emphasized the fact that the set of the attractors generated by
GIFSs is larger than the one generated by iterated function systems. Another
related topics can be found in [4], [6], [8], [10], [11], [12], [13] and
[16].

Moreover, in [3], Strobin and his collaborators provided four algorithms
which generate images of attractors of GIFSs, one of them being the
deterministic algorithm for GIFSs (a counterpart of the classical
deterministic algorithm for iterated function systems). Note that in [3] one
can also find a list of papers dealing with algorithms generating images of
the attractors of iterated function systems.

In this paper we present another algorithm (called the grid algorithm)
allowing to generate images of the attractors of GIFSs on finite dimensional
spaces and we compare it with the deterministic algorithm for GIFSs. The
deterministic algorithm for GIFSs consists in choosing a finite set of
points and applying to this set each of the constitutive functions of the
system obtaining in this way a new finite set of points. To each of these
new points we apply again each of the constitutive functions of the system.
Continuing the procedure described above we approach the attractor. The main
idea of the grid algorithm is to simplify the deterministic algorithm by
dividing, at each step, the space that we are working with in small pieces
and to choose for each such piece just one point.

\bigskip

\textbf{II. Preliminaries}

\bigskip

Given a metric space $(X,d)$, we adopt the following notation: 
\begin{equation*}
P_{cp}(X)\overset{not}{=}\{K\subseteq X\mid K\text{ }\ \text{is non-empty
and compact}\}\text{.}
\end{equation*}

For $K_{1},K_{2}\in P_{cp}(X)$, we consider 
\begin{equation*}
d(K_{1},K_{2})\overset{def}{=}\underset{x\in K_{1}}{\sup }d(x,K_{2})\text{,}
\end{equation*}
where $d(x,K_{2})\overset{def}{=}\underset{y\in K_{2}}{\inf }d(x,y)$.

The function $h:P_{cp}(X)\times P_{cp}(X)\rightarrow \lbrack 0,\infty )$
given by 
\begin{equation*}
h(K_{1},K_{2})=\max \{d(K_{1},K_{2}),d(K_{2},K_{1})\}\text{,}
\end{equation*}
for every $K_{1},K_{2}\in P_{cp}(X)$, turns out to be a metric which is
called the Hausdorff-Pompeiu metric.

If $(X,d)$ is complete, then $(P_{cp}(X),h)$ is complete.

\bigskip

Given a metric space $(X,d)$ and $p\in \mathbb{N}^{\ast }$, by $X^{p}$ we
denote the Cartesian product of $p$ copies of $X$. We endow $X^{p}$ with the
maximum metric $d_{\max }$ defined by 
\begin{equation*}
d_{\max }((x_{1},...,x_{p}),(y_{1},...,y_{p}))=\max
\{d(x_{1},y_{1}),...,d(x_{p},y_{p})\}\text{,}
\end{equation*}%
for all $(x_{1},...,x_{p}),(y_{1},...,y_{p})\in X^{p}$.

\bigskip

\textbf{Definition 2.1}. \textit{A generalized iterated function system (of
order }$p$\textit{) is a pair }$\mathcal{F}=((X,d),(f_{i})_{i\in
\{1,2,...,L\}})$\textit{, where }$(X,d)$ \textit{is a metric space, }$p,L\in 
\mathbb{N}^{\ast }$\textit{\ and }$f_{i}:X^{p}\rightarrow X$\textit{\ is
contraction for each }$i\in \{1,...,L\}$\textit{. The function }$F
_{\mathcal{F}}:(P_{cp}(X))^{p}\rightarrow P_{cp}(X)$, \textit{described by} 
\begin{equation*}
F_{\mathcal{F}}(K_{1},...,K_{p})=\underset{i\in \{1,...,L\}}{\cup }%
f_{i}(K_{1}\times ...\times K_{p})\text{,}
\end{equation*}%
\textit{for all} $K_{1},...,K_{p}\in P_{cp}(X)$, \textit{is called the
fractal operator associated to} $\mathcal{F}$.

\bigskip

We shall use the abbreviation GIFS for a generalized iterated function
system.

\bigskip

\textbf{Theorem 2.2} (see Theorem 3.9 from [9])\textbf{.} \textit{Given a
complete metric space }$(X,d)$ \textit{and a GIFS} $\mathcal{F=}%
((X,d),(f_{i})_{i\in \{1,...,L\}})$ \textit{of order} $p$\textit{, there
exists a unique} $A_{\mathcal{F}}\in P_{cp}(X)$ \textit{such that} 
\begin{equation*}
F_{\mathcal{F}}(A_{\mathcal{F}},...,A_{\mathcal{F}})=A_{\mathcal{F}%
}\text{\textit{.}}
\end{equation*}

\textit{In addition, for every} $K_{1},...,K_{p}\in P_{cp}(X)$\textit{, the
sequence} $(K_{n})_{n}$ \textit{defined by} 
\begin{equation*}
K_{n+p}=F_{\mathcal{F}}(K_{n},...,K_{n+p-1})\text{,}
\end{equation*}%
\textit{for every} $n\in \mathbb{N}^{\ast }$, \textit{converges, with
respect to the Hausdorff-Pompeiu metric, to} $A_{\mathcal{F}}$.

\bigskip

\textbf{Definition 2.3}. \textit{In the framework of the above theorem, the
set} $A_{\mathcal{F}}$ \textit{is called the fractal generated by} $\mathcal{%
F}$\textit{.}

\bigskip

\textbf{Remark 2.4 }(see Remark 12 from [3]). \textit{In the framework of
the above definition, the function} $\mathcal{G}_{\mathcal{F}%
}:P_{cp}(X)\rightarrow P_{cp}(X)$, \textit{described by}%
\begin{equation*}
\mathcal{G}_{\mathcal{F}}(K)=F_{\mathcal{F}}(K,...,K)=\underset{%
i\in \{1,...,L\}}{\cup }f_{i}(K\times...\times K)\text{,}
\end{equation*}%
\textit{for all} $K\in P_{cp}(X)$,\textit{\ is a contraction on the complete
metric space} $(P_{cp}(X),h)$\textit{\ since it has the Lipschitz constant
less of equal to} $\max \{lip(f_{1}),...,lip(f_{L})\}<1$\textit{.}

\bigskip

\textbf{III.} \textbf{The presentation of the algorithms}

\bigskip

For $(x_{1},....,x_{M})\in \mathbb{R}^{M}$, we shall use the following
notation:%
\begin{equation*}
\lbrack (x_{1},....,x_{M})]=([x_{1}],....,[x_{M}])\text{,}
\end{equation*}%
where $[x]$ designates the greatest integer less than or equal to the real
number $x$.

\medskip

In the sequel, without loss of generality, 
\begin{equation*}
\mathcal{F}=(([0,D]^{M},d),\{f_{1},...,f_{L}\})\text{,}
\end{equation*}%
where $L,M\in \mathbb{N}$ and $d$ is the euclidean distance in $\mathbb{R}%
^{M}$, will be a generalized iterated function system of order $p\geq 2$ (so 
$f_{i}:([0,D]^{M})^{p}\rightarrow \lbrack 0,D]^{M}$ for every $i\in
\{1,...,L\}$).

\medskip

We shall use the following notation:

$\bullet $ $\max \{lip(f_{1}),...,lip(f_{L})\}\overset{not}{=}C<1$

$\bullet $ $\beta =pM$.

\medskip

We also consider the following functions:

$\bullet $ $ F_{\mathcal{F}}:(P_{cp}([0,D]^{M}))^{p}\rightarrow
P_{cp}([0,D]^{M})$ described by%
\begin{equation*}
 F_{\mathcal{F}}(K_{1},...,K_{p})=f_{1}(K_{1}\times...\times K_{p})\cup
...\cup f_{L}(K_{1}\times...\times K_{p})\text{,}
\end{equation*}%
for all $K_{1},...,K_{p}\in P_{cp}([0,D]^{M})$

$\bullet $ $\mathcal{G}_{\mathcal{F}}:P_{cp}([0,D]^{M})\rightarrow
P_{cp}([0,D]^{M})$ described by%
\begin{equation*}
\mathcal{G}_{\mathcal{F}}(K)= F_{\mathcal{F}}(K,...,K)\text{,}
\end{equation*}%
for every $K\in P_{cp}([0,D]^{M})$.

$\bullet $ $(n_{k})_{k\in \mathbb{N}^{\ast }}$ a sequence of natural numbers.

\medskip

For a finite set $K_{0}\in P_{cp}([0,D]^{M})$, we shall use the following
notations:

$\bullet $ 
\begin{equation*}
A_{k}\overset{not}{=}\mathcal{G}_{\mathcal{F}}^{[k]}(K_{0})\text{,}
\end{equation*}%
where $k\in \mathbb{N}$

$\bullet $%
\begin{equation*}
\overset{\sim }{A_{k}}\overset{not}{=}\{\frac{D}{n_{k}}[\frac{n_{k}}{D}%
f_{l}(u_{1},...,u_{p})]\mid u_{1},...,u_{p}\in \overset{\sim }{A}_{k-1}\text{%
, }l\in \{1,...,L\}\}\text{,}
\end{equation*}%
where $k\in \mathbb{N}^{\ast }$ and $\overset{\sim }{A_{0}}=K_{0}$

$\bullet $%
\begin{equation*}
\frac{D\sqrt{M}}{n_{k}}\overset{not}{=}\varepsilon _{k}\text{,}
\end{equation*}%
where $k\in \mathbb{N}$.

\bigskip

Let us recall the pseudocode for the deterministic algorithm for a GIFS (see
[3])

\bigskip

\textbf{Pseudocode for the deterministic algorithm for a GIFS}

\bigskip

Read initially defined objects: constants: $L,M$, finite set, $m$ natural
number: $K_{0}\in P_{cp}([0,D]^{M})$, mappings: $f_{1},...,f_{L}$,
variables: $k,D_{0}$.

Initial values: $D_{0}:=K_{0}$.

\qquad \qquad \qquad For $k$ from $1$ to $m-1$

\qquad \qquad \qquad \qquad $D_{1}:=\mathcal{G}_{\mathcal{F}}(D_{0})$

\qquad \qquad \qquad \qquad $D_{0}:=D_{1}$.

Print $D_{m}$.

\bigskip

Now let us present the pseudocode for our new algorithm.

\bigskip

\textbf{Pseudocode for the grid algorithm for a GIFS}

\bigskip

Read initially defined objects: constant: $L,$ $M$, finite set, $m$ natural
number: $K_{0}\in P_{cp}([0,D]^{M})$, mappings: $f_{1},...,f_{L}$, sequence: 
$(n_{k})_{k}$, variables: $k,D_{0}$.

Initial values: $D_{0}:=K_{0}$.

\qquad \qquad For $k$ from $1$ to $m-1$

\qquad \qquad \qquad $D_{1}:=\{\frac{D}{n_{k}}[\frac{n_{k}}{D}%
f_{l}(u_{1},...,u_{p})]\mid u_{1},...,u_{p}\in D_{0},l\in \{1,...,L\}\}$

\qquad \qquad \qquad $D_{0}:=D_{1}$.

Print $D_{m}$.

\newpage

\textbf{IV.} \textbf{The complexity of the algorithms}

\bigskip

By $x_{k}$ we denote the number of points computed at the step $k$ of the
deterministic algorithm and by $y_{k}$ the number of points computed up to
the step $k$ of the grid algorithm.

\bigskip

\textbf{A}. \textit{The case of the deterministic algorithm}

\bigskip

We have $x_{k+1}\leq L(x_{k})^{p}$, so, with the notation $z_{k}\overset{not}%
{=}\ln x_{k}$, we obtain $z_{k+1}\leq \ln L+pz_{k}$ for every $k\in \mathbb{N%
}$. Therefore $z_{k}\leq \frac{p^{k}-1}{p-1}\ln L+p^{k}z_{0}$, i.e.%
\begin{equation}
x_{k}\leq \frac{1}{L^{\frac{1}{p-1}}}(x_{0}L^{\frac{1}{p-1}})^{p^{k}}\text{,}
\tag{1}
\end{equation}%
for every $k\in \mathbb{N}$.

Note that, according to Remark 2.4, we have $h(A_{k},A_{\mathcal{F}})\leq 
\frac{h(A_{0},A_{1})}{1-C}C^{k}$ for every $k\in \mathbb{N}$. Therefore, in
order to be sure that $A_{k}$ approximates the attractor $A_{\mathcal{F}}$
with accuracy $\varepsilon \frac{h(A_{0},A_{1})}{1-C}$, we need $k>\log
_{C^{-1}}(\varepsilon ^{-1})$. Hence, based on $(1)$, the quantity $\frac{1}{%
L^{\frac{1}{p-1}}}(x_{0}L^{\frac{1}{p-1}})^{p^{\log _{C^{-1}}(\varepsilon
^{-1})}}=\frac{1}{L^{\frac{1}{p-1}}}(x_{0}L^{\frac{1}{p-1}})^{(\frac{1}{%
\varepsilon })^{\log _{C^{-1}}(p)}}$ describes the number of points that we
have to compute in order to be sure that $A_{k}$ is an approximation of $A_{%
\mathcal{F}}$ with an error less than $\varepsilon \frac{h(A_{0},A_{1})}{1-C}
$.

\bigskip

\textbf{Conclusion}: \textit{The complexity of the deterministic algorithm
is described by the function }$\mathcal{C}_{c}:(0,\infty )\rightarrow 
\mathbb{R}$\textit{\ given by}%
\begin{equation*}
\mathcal{C}_{c}(\varepsilon )=(x_{0}L^{\frac{1}{p-1}})^{(\frac{1}{%
\varepsilon })^{\log _{\frac{1}{C}}(p)}}\text{,}
\end{equation*}%
\textit{for every }$\varepsilon >0$\textit{.}

\bigskip

\textbf{B}. \textit{The case of the grid algorithm}

\bigskip

\textbf{Remark 4.1}. Since $y_{k+1}\leq L(n_{k})^{\beta }$ for every $k\in 
\mathbb{N}^{\ast }$, \textit{up to the step }$N$\textit{, we have to compute}
$L\overset{N}{\underset{k=1}{\sum }}(n_{k})^{\beta }=L(D\sqrt{M})^{\beta }%
\overset{N}{\underset{k=1}{\sum }}(\frac{1}{\varepsilon _{k}})^{\beta }$ 
\textit{points.}

\bigskip

\textbf{Remark 4.2}. \textit{We have}%
\begin{equation*}
h(\overset{\sim }{A_{k}},\mathcal{G}_{\mathcal{F}}(\overset{\sim }{A}%
_{k-1}))\leq \varepsilon _{k}\text{,}
\end{equation*}%
\textit{for every} $k\in \mathbb{N}^{\ast }$.

\bigskip

\textbf{Remark 4.3}. \textit{We have} 
\begin{equation*}
h(\overset{\sim }{A_{0}},A_{\mathcal{F}})\leq diam([0,D]^{M})=D\sqrt{M}\text{%
.}
\end{equation*}

\bigskip

As the inequality%
\begin{equation*}
h(\overset{\sim }{A_{k}},A_{\mathcal{F}})\leq h(\overset{\sim }{A_{k}},%
\mathcal{G}_{\mathcal{F}}(\overset{\sim }{A}_{k-1}))+h(\mathcal{G}_{\mathcal{%
F}}(\overset{\sim }{A}_{k-1}),\mathcal{G}_{\mathcal{F}}(A_{\mathcal{F}}))\leq
\end{equation*}%
\begin{equation*}
\overset{\text{Remarks 2.4 and 4.2}}{\leq }\varepsilon _{k}+Ch(\overset{\sim 
}{A}_{k-1},A_{\mathcal{F}})\text{,}
\end{equation*}%
is valid for every $k\in \mathbb{N}^{\ast }$, we get%
\begin{equation*}
h(\overset{\sim }{A_{k}},A_{\mathcal{F}})\leq \varepsilon _{k}+C\varepsilon
_{k-1}+C^{2}\varepsilon _{k-2}+...+C^{k-2}\varepsilon
_{2}+C^{k-1}\varepsilon _{1}+C^{k}h(\overset{\sim }{A_{0}},A_{\mathcal{F}})%
\text{,}
\end{equation*}%
so, taking into account Remark 4.3, we obtain%
\begin{equation*}
h(\overset{\sim }{A_{k}},A_{\mathcal{F}})\leq \varepsilon _{k}+C\varepsilon
_{k-1}+C^{2}\varepsilon _{k-2}+...+C^{k-2}\varepsilon
_{2}+C^{k-1}\varepsilon _{1}+C^{k}D\sqrt{M}\text{,}
\end{equation*}%
for every $k\in \mathbb{N}^{\ast }$. Consequently, we arrive to the
following problem: given a fixed natural number $N$ and $\varepsilon >0$
such that $\frac{\varepsilon }{C^{N}}-D\sqrt{M}>0$, find the minimum of the
function $f:[0,\infty )^{N}\rightarrow \lbrack 0,\infty )$, given by%
\begin{equation*}
f(\varepsilon _{1},...,\varepsilon _{N})=\overset{N}{\underset{k=1}{\sum }}(%
\frac{1}{\varepsilon _{k}})^{\beta }\text{,}
\end{equation*}%
for every $\varepsilon _{1},...,\varepsilon _{N}\in \lbrack 0,\infty )$,
with the constraint 
\begin{equation*}
\varepsilon _{N}+C\varepsilon _{N-1}+C^{2}\varepsilon
_{N-2}+...+C^{N-2}\varepsilon _{2}+C^{N-1}\varepsilon _{1}+C^{N}D\sqrt{M}%
=\varepsilon \text{.}
\end{equation*}

We adopt the following notations:

$\bullet $ $t\overset{not}{=}C^{-\frac{\beta }{\beta +1}N}-1$

$\bullet $ $K_{1}\overset{not}{=}\frac{C^{\frac{1}{\beta +1}}-C}{C^{\frac{1}{%
\beta +1}}}=1-C^{\frac{\beta }{\beta +1}}$

$\bullet $ $K_{2}\overset{not}{=}K_{1}^{-\beta -1}$

$\bullet $ $K_{3}\overset{not}{=}K_{2}\varepsilon ^{-\beta }$

$\bullet $ $a\overset{not}{=}\frac{D\sqrt{M}}{\varepsilon }$

$\bullet $ $y\overset{not}{=}\frac{1}{C^{N}}$.

Since we are going to use the method of Lagrange multipliers, we consider
the function $F=f+\lambda g$, where $\lambda \in \mathbb{R}$ and the
function $g:[0,\infty )^{N}\rightarrow \lbrack 0,\infty )$ is given by%
\begin{equation*}
g(\varepsilon _{1},...,\varepsilon _{N})=\varepsilon _{N}+C\varepsilon
_{N-1}+...+C^{N-2}\varepsilon _{2}+C^{N-1}\varepsilon _{1}+C^{N}D\sqrt{M}%
-\varepsilon \text{,}
\end{equation*}%
for every $\varepsilon _{1},...,\varepsilon _{N}\in \lbrack 0,\infty )$. The
equation $\frac{\partial F}{\partial \varepsilon _{k}}=0$, i.e. $-\beta
(\varepsilon _{k})^{-\beta -1}+\lambda C^{N-k}=0$, has the solution%
\begin{equation}
\varepsilon _{k}^{0}=k_{N}C^{\frac{k}{\beta +1}}\text{,}  \tag{1}
\end{equation}%
for every $k$, where $k_{N}=\frac{1}{C^{\frac{N}{\beta +1}}}(\frac{\beta }{%
\lambda })^{\frac{1}{\beta +1}}$. As $g(\varepsilon _{1}^{0},...,\varepsilon
_{N}^{0})=0$, we get 
\begin{equation*}
k_{N}(C^{\frac{N}{\beta +1}}+C^{1+\frac{N-1}{\beta +1}}+C^{2+\frac{N-2}{%
\beta +1}}+...+C^{N-2+\frac{2}{\beta +1}}+C^{N-1+\frac{1}{\beta +1}%
})=\varepsilon -C^{N}D\sqrt{M}\text{,}
\end{equation*}%
i.e.%
\begin{equation*}
k_{N}C^{\frac{N}{\beta +1}}(1+C^{\frac{\beta }{\beta +1}}+C^{^{2\frac{\beta 
}{\beta +1}}}+...+C^{^{^{(N-2)\frac{\beta }{\beta +1}}}}+C^{(N-1)\frac{\beta 
}{\beta +1}})=\varepsilon -C^{N}D\sqrt{M}\text{,}
\end{equation*}%
so $k_{N}C^{\frac{N}{\beta +1}}\frac{(C^{\frac{\beta }{\beta +1}})^{N}-1}{C^{%
\frac{\beta }{\beta +1}}-1}=\varepsilon -C^{N}D\sqrt{M}$, which implies $%
k_{N}\frac{C^{N}-C^{\frac{N}{\beta +1}}}{C-C^{\frac{1}{\beta +1}}}=\frac{%
\varepsilon -C^{N}D\sqrt{M}}{C^{\frac{1}{\beta +1}}}$. The last equality
takes the form $k_{N}\frac{C^{-\frac{\beta }{\beta +1}N}-1}{C^{\frac{1}{%
\beta +1}}-C}=\frac{\frac{\varepsilon }{C^{N}}-D\sqrt{M}}{C^{\frac{1}{\beta
+1}}}$. Thus we obtain%
\begin{equation}
k_{N}=\frac{K_{1}}{t}(\frac{\varepsilon }{C^{N}}-D\sqrt{M})\text{.}  \tag{2}
\end{equation}%
We have%
\begin{equation*}
f(\varepsilon _{1}^{0},...,\varepsilon _{N}^{0})=\overset{N}{\underset{k=1}{%
\sum }}(\varepsilon _{k}^{0})^{-\beta }\overset{(1)}{=}(k_{N})^{-\beta }%
\overset{N}{\underset{k=1}{\sum }}C^{-\frac{\beta }{\beta +1}%
k}=(k_{N})^{-\beta }C^{-\frac{\beta }{\beta +1}}\frac{(C^{-\frac{\beta }{%
\beta +1}})^{N}-1}{C^{-\frac{\beta }{\beta +1}}-1}=
\end{equation*}%
\begin{equation*}
\overset{(2)}{=}t^{\beta }K_{1}^{-\beta }(\frac{\varepsilon }{C^{N}}-D\sqrt{M%
})^{-\beta }\frac{t}{1-C^{\frac{\beta }{\beta +1}}}=t^{\beta +1}(\frac{%
\varepsilon }{C^{N}}-D\sqrt{M})^{-\beta }\frac{K_{1}^{-\beta }}{K_{1}}=
\end{equation*}%
\begin{equation*}
=t^{\beta +1}(\frac{\varepsilon }{C^{N}}-D\sqrt{M})^{-\beta }K_{1}^{-\beta
-1}=K_{2}(\frac{\varepsilon }{C^{N}}-D\sqrt{M})^{-\beta }(C^{-\frac{\beta }{%
\beta +1}N}-1)^{\beta +1}\text{.}
\end{equation*}%
Therefore, the last equality can be written as%
\begin{equation}
f(\varepsilon _{1}^{0},...,\varepsilon _{N}^{0})=K_{3}(y^{\frac{\beta }{%
\beta +1}}-1)^{\beta +1}(y-a)^{-\beta }\text{.}  \tag{3}
\end{equation}

As the right hand side of $(3)$ gives us the optimal number of points that
we have to compute, after $N$ steps, in order to approximate $A_{\mathcal{F}%
} $ by $\overset{\sim }{A_{k}}$ with an error not greater than $\varepsilon $%
, we need to find the minimum value of the function $h:(a,\infty
)\rightarrow \mathbb{R}$ given by 
\begin{equation*}
h(y)=K_{3}(y^{\frac{\beta }{\beta +1}}-1)^{\beta +1}(y-a)^{-\beta }\text{,}
\end{equation*}%
for every $y\in (a,\infty )$. One can easily see that%
\begin{equation*}
h^{^{\prime }}(y)=K_{3}\beta (y^{\frac{\beta }{\beta +1}}-1)^{\beta
}(y-a)^{-\beta -1}(1-ay^{-\frac{1}{\beta +1}})\text{,}
\end{equation*}%
for every $y\in (a,\infty )$. As we can suppose that $a>1$ (since we are
interested in the case when $\varepsilon $ is small), $\underset{%
y\rightarrow \infty }{\lim }h(y)=K_{3}$ and $\underset{\underset{y>a}{%
y\rightarrow a}}{\lim }h(y)=\infty $, we conclude that $h$ attains its
minimum at $a^{\beta +1}$ and the value of the minimum is 
\begin{equation*}
h(a^{\beta +1})=K_{3}(a^{\beta }-1)^{\beta +1}(a^{\beta +1}-a)^{-\beta }=
\end{equation*}%
\begin{equation*}
=K_{3}\frac{a^{\beta }-1}{a^{\beta }}=\varepsilon ^{-\beta }(1-C^{\frac{%
\beta }{\beta +1}})^{-\beta -1}\frac{(\frac{D\sqrt{M}}{\varepsilon })^{\beta
}-1}{(\frac{D\sqrt{M}}{\varepsilon })^{\beta }}\text{,}
\end{equation*}%
so%
\begin{equation*}
\underset{\underset{\varepsilon >0}{\varepsilon \rightarrow 0}}{\lim }\frac{%
h(a^{\beta +1})}{(1-C^{\frac{\beta }{\beta +1}})^{-\beta -1}(\frac{1}{%
\varepsilon })^{\beta }}=1\text{.}
\end{equation*}

\bigskip

\textbf{Conclusion}: \textit{The complexity of the grid algorithm is
described by the function }$\mathcal{C}_{g}:(0,\infty )\rightarrow \mathbb{R}
$\textit{\ given by}%
\begin{equation*}
\mathcal{C}_{g}(\varepsilon )=(1-C^{\frac{\beta }{\beta +1}})^{-\beta -1}(%
\frac{1}{\varepsilon })^{pM}\text{,}
\end{equation*}%
\textit{for every }$\varepsilon >0$\textit{.}

\bigskip

In the final of this section we mention that (in order to avoid very
complicated computations) we did not pay attention to the fact that the best
values of $n_{k}$ and $N$ that we obtained (namely $\frac{D\sqrt{M}}{%
\varepsilon _{k}^{0}}$ and $(\beta +1)\frac{\ln (\frac{\varepsilon }{D\sqrt{M%
}})}{\ln (C)}$) may not be integers. In reality we could work with $n_{k}=[%
\frac{D\sqrt{M}}{\varepsilon _{k}^{0}}]+1$ and $N=[(\beta +1)\frac{\ln (%
\frac{\varepsilon }{D\sqrt{M}})}{\ln (C)}]+1$ without a significant change.

\bigskip

\textbf{V. Final remarks}

\bigskip

\textbf{Remark 5.1}. We have 
\begin{equation*}
\underset{\underset{\varepsilon >0}{\varepsilon \rightarrow 0}}{\lim }\frac{%
\mathcal{C}_{g}(\varepsilon )}{\mathcal{C}_{c}(\varepsilon )}=\underset{%
\underset{\varepsilon >0}{\varepsilon \rightarrow 0}}{\lim }\frac{(1-C^{%
\frac{\beta }{\beta +1}})^{-\beta -1}(\frac{1}{\varepsilon })^{pM}}{(x_{0}L^{%
\frac{1}{p-1}})^{(\frac{1}{\varepsilon })^{\log _{\frac{1}{C}}(p)}}}=0\text{,%
}
\end{equation*}%
\textit{so the grid algorithm is more efficient than the deterministic
algorithm.}

\bigskip

\textbf{Remark 5.2}. As $\left\vert u-[u+\frac{1}{2}]\right\vert \leq \frac{1%
}{2}$, \textit{we can improve our grid algorithm} (which is based on the
inequality $\left\vert u-[u]\right\vert <1$)\textit{\ in the following way:}

\bigskip

\textbf{Pseudocode for the improved grid algorithm for GIFS}

\bigskip

Read initially defined objects: constant: $L,$ $M$, finite set, $m$ natural
number: $K_{0}\in P_{cp}([0,D]^{M})$, mappings: $f_{1},...,f_{L}$, sequence: 
$(n_{k})_{k}$, variables: $k,D_{0}$.

Initial values: $D_{0}:=K_{0}$.

\qquad \qquad For $k$ from $1$ to $m-1$

\qquad \qquad \qquad $D_{1}:=\{\frac{D}{n_{k}}[\frac{n_{k}}{D}%
f_{l}(u_{1},...,u_{p})+\frac{1}{2}]\mid u_{1},...,u_{p}\in D_{0},l\in
\{1,...,L\}\}$

\qquad \qquad \qquad $D_{0}:=D_{1}$.

Print $D_{m}$.

\bigskip

\textbf{Remark 5.3}. On the one hand, repeating the arguments used in IV, A,
for the case of an iterated function system (i.e. $p=1$), we obtain that the
complexity of the corresponding algorithm is\textit{\ }described by the
function\textit{\ }$\mathcal{C}:(0,\infty )\rightarrow \mathbb{R}$\textit{\ }%
given by%
\begin{equation*}
\mathcal{C}(\varepsilon )=(\frac{1}{\varepsilon })^{\frac{\ln L}{\ln \frac{1%
}{C}}}\text{,}
\end{equation*}%
for every\textit{\ }$\varepsilon >0$\textit{, }so $C$ is involved at the
exponent of $\frac{1}{\varepsilon }$.\textit{\ }We stress upon the fact that
since $\underset{\underset{C<1}{C\rightarrow 1}}{\lim }\frac{1}{\ln \frac{1}{%
C}}=\infty $, the closer is $C$ to $1$, the bigger is the number of points
that we have to compute in order to approximate the attractor with an error
less that $\varepsilon $.

On the other hand, in the rule that gives $\mathcal{C}_{g}(\varepsilon )$
the constant $C$ is involved only in the coefficient $(1-C^{\frac{\beta }{%
\beta +1}})^{-\beta -1}$.

Moreover, note that $\underset{\underset{C<1}{C\rightarrow 1}}{\lim }\frac{%
\mathcal{C}_{g}(\varepsilon )}{\mathcal{C}_{c}(\varepsilon )}=0$ for each $%
\varepsilon \in (0,1)$.

\bigskip

\textbf{VI. Examples}

\bigskip 

In section IV we compared the algorithms with respect to a fixed preassigned
error. In this section our goal is to get an optimal image (with respect to
the computer that we worked with) for three examples. For this reason we
chose a version of the grid algorithm for which $n_{k}=k^{2}$, the error
being less than $D\sqrt{M}\left(\frac{1}{n^{2}}+C\frac{1}{(n-1)^{2}}+C^{2}\frac{1}{(n-2)^{2}%
}+...+C^{n}\right)$, where $n$ is the number of steps and $C$ is the contraction
constant of the system.

\bigskip

\textbf{A}.

Consider the GIFS $\mathcal{F}=(([0,1]^{2},d),\{f_{1},f_{2},f_{3}\})$,
where, for $x=(x_{1},y_{1})$ and $y=(x_{2},y_{2})$, we have 
\begin{equation*}
f_{1}(x,y)=(0.2x_{1}+0.2y_{2};0.2x_{2}+0.1y_{2})
\end{equation*}%
\begin{equation*}
f_{2}(x,y)=(0.15x_{1}+0.07x_{2}+0.4;0.15y_{1}+0.07y_{2})\text{.}
\end{equation*}%
and%
\begin{equation*}
f_{3}(x,y)=(0.15y_{1}+0.07x_{2};0.15x_{1}+0.07y_{2}+0.04)\text{.}
\end{equation*}

Using the deterministic algorithm we get the image indicated in figure 1
and using the grid algorithm we get the image in figure 2.

The deterministic algorithm run 4 steps in 10 seconds, while the grid
algorithm 8 steps in less that 10 seconds.

Figure 1 has approximately $3^{2^{4}}=43046721$ points, while figure 2
comprises around 20000 points.

\bigskip

\textbf{Remark 6.1}. \textit{If we allow the deterministic algorithm to run
5 steps it needs about 90 minutes and we get a very similar image with the
one in figure 2.}

\bigskip

\begin{center}

\includegraphics{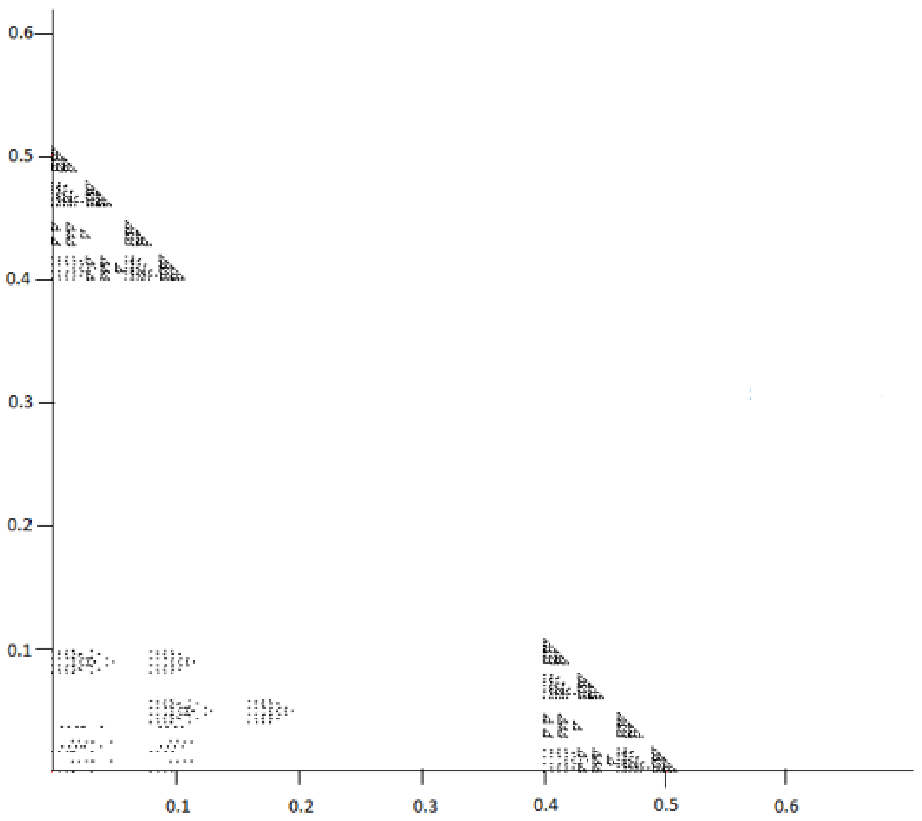}
\\
Figure 1
\end{center}

\begin{center}
\includegraphics[width=9.5cm]
{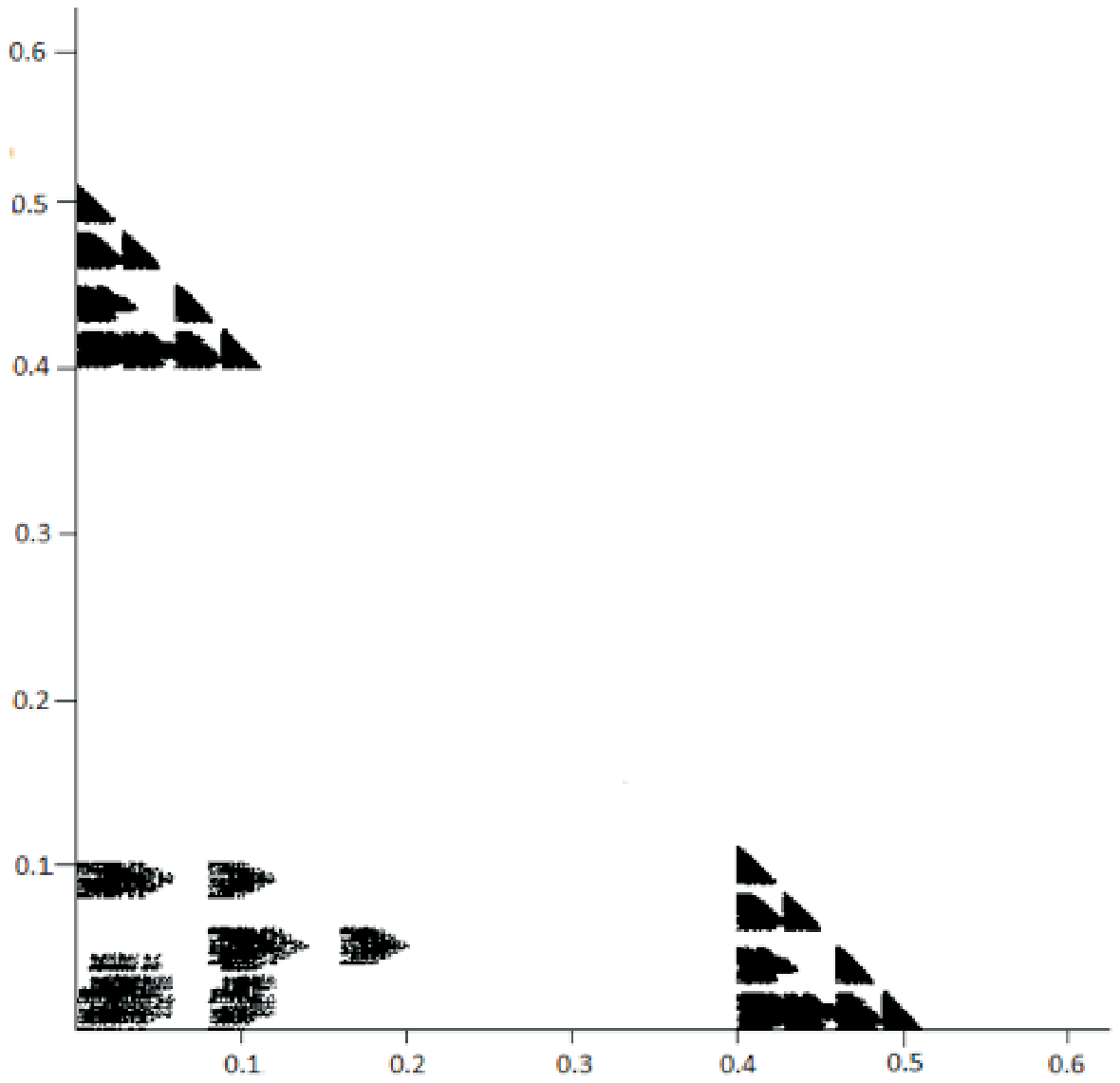}
\\
 Figure 2
\end{center}

\bigskip

\textbf{B}.

Consider the GIFS $\mathcal{F}=(([0,1]^{2},d),\{f_{1},f_{2}\})$, where, for $%
x=(x_{1},y_{1})$ and $y=(x_{2},y_{2})$, we have 
\begin{equation*}
f_{1}(x,y)=(0.1x_{1}+0.16y_{1}-0.01x_{2}+0.3y_{2};-0.05y_{1}+0.15x_{2}+0.15y_{2})
\end{equation*}%
and%
\begin{equation*}
f_{2}(x,y)=(0.09x_{1}-0.1y_{1}-0.15x_{2}+0.14y_{2}+0.4;0.14x_{1}+0.14y_{1}+0.14x_{2}+0.04)%
\text{.}
\end{equation*}

The deterministic algorithm yields the image in figure 3 and the grid
algorithm produces the image in figure 4.

The deterministic algorithm needed 20 seconds to run 5 steps, while the grid
algorithm needed about 10 minutes to run 14 steps.

Figure 3 consists of about $2^{2^{5}}=4294967296$ points, while figure 4 is built up using around $2000000$ p oints.

\bigskip

\textbf{Remark 6.2}. \textit{With the aid of the computer that we utilized,
the deterministic algorithm would need 42 days to run 6 steps.}

\bigskip

\begin{center}
\includegraphics[width=8.5cm]{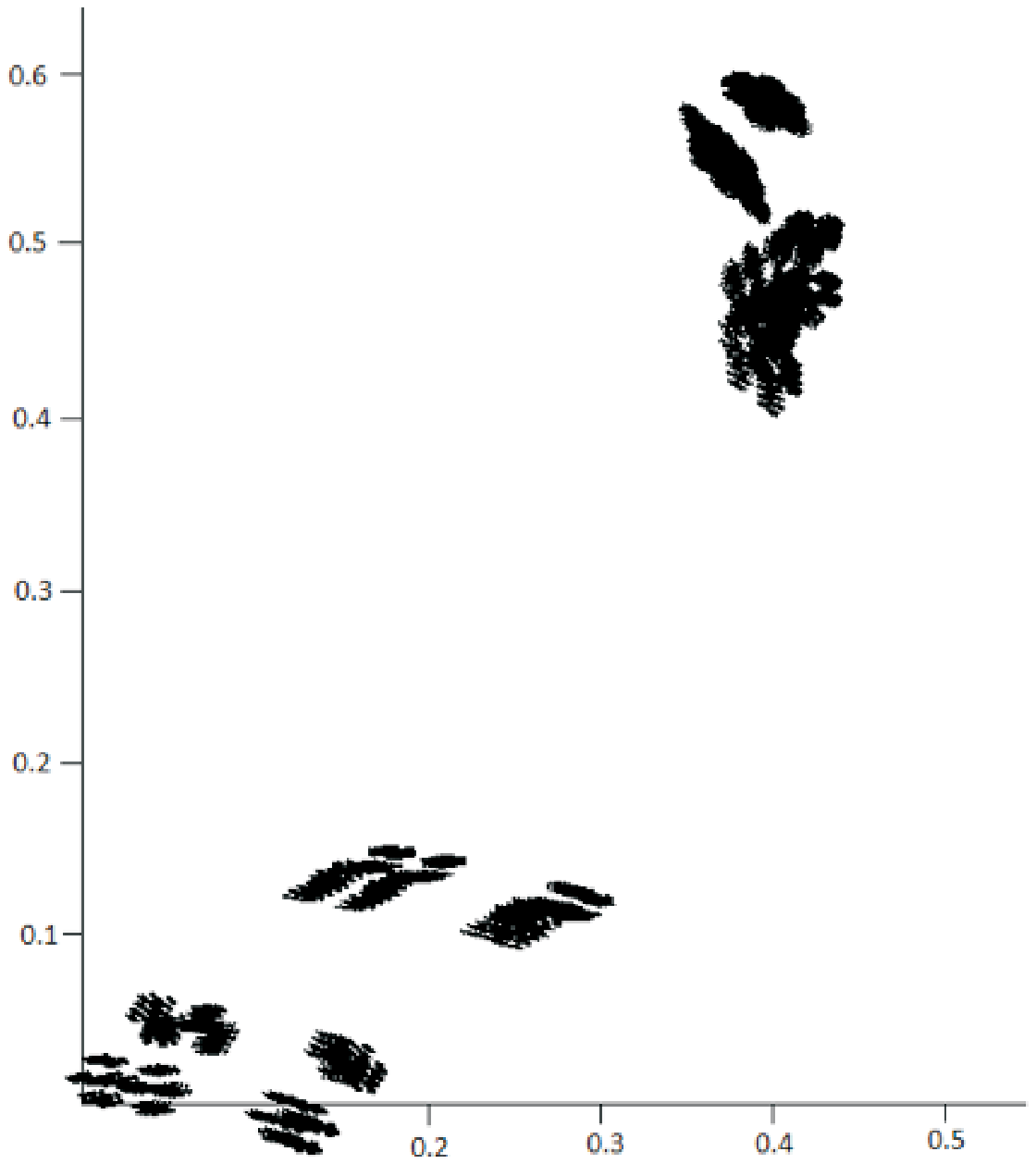}\\
Figure 3
\end{center}

\begin{center}
\includegraphics[width=8cm]{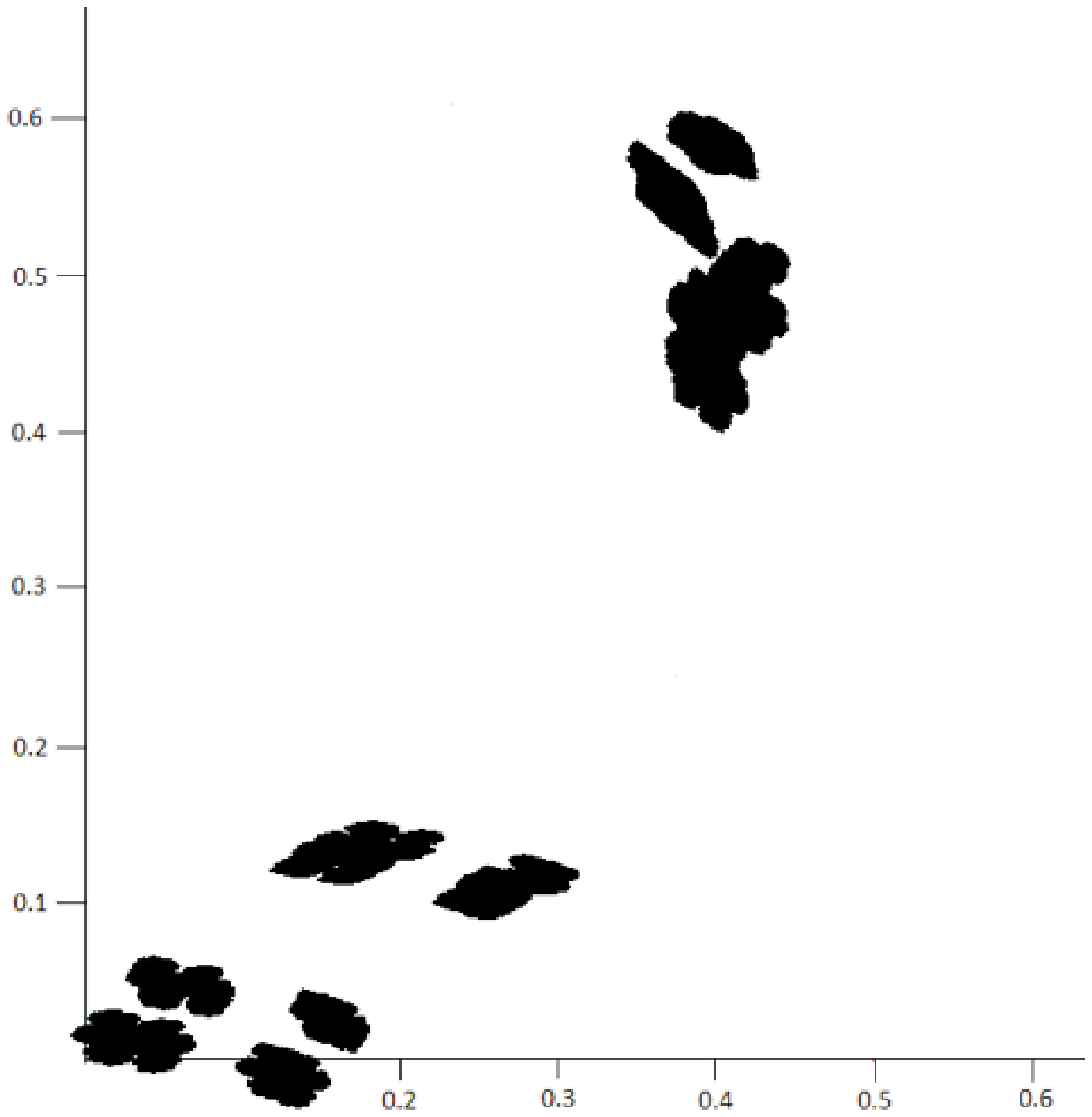}\\
Figure 4
\end{center}

\bigskip

\textbf{C}.

Consider the GIFS $\mathcal{F}=(([0,1]^{2},d),\{f_{1},f_{2}\})$, where, for $%
x=(x_{1},y_{1})$ and $y=(x_{2},y_{2})$, we have 
\begin{equation*}
f_{1}(x,y)=(0.5x_{1}-0.5y_{1}+0.001x_{2}+0.45;0.5x_{1}+0.5y_{1}+0.001y_{2}-0.05)
\end{equation*}%
and%
\begin{equation*}
f_{2}(x,y)=(0.2x_{1}+0.01x_{2}+0.14y_{2}+0.147;0.2y_{1}+0.01y_{2}+0.105)%
\text{.}
\end{equation*}

The image in figure 5 indicates what we get running the deterministic
algorithm and the image in figure 6 what we obtain using the grid algorithm.

Both algorithms ran about 2 minutes, the deterministic one running 5 steps,
while the grid one 22 steps.

Figure 5 consists of about $2^{2^{5}}=4294967296$ points, while figure 6
is made up of circa $217800$ points.

\medskip

\textbf{Remark 6.3}. \textit{Even though the number of points making up
figure 5 is considerably bigger that the number of points building up
figure 6, one can observe that the grid algorithm produces a much better
approximation of the attractor.}

\bigskip

\begin{center}
\includegraphics{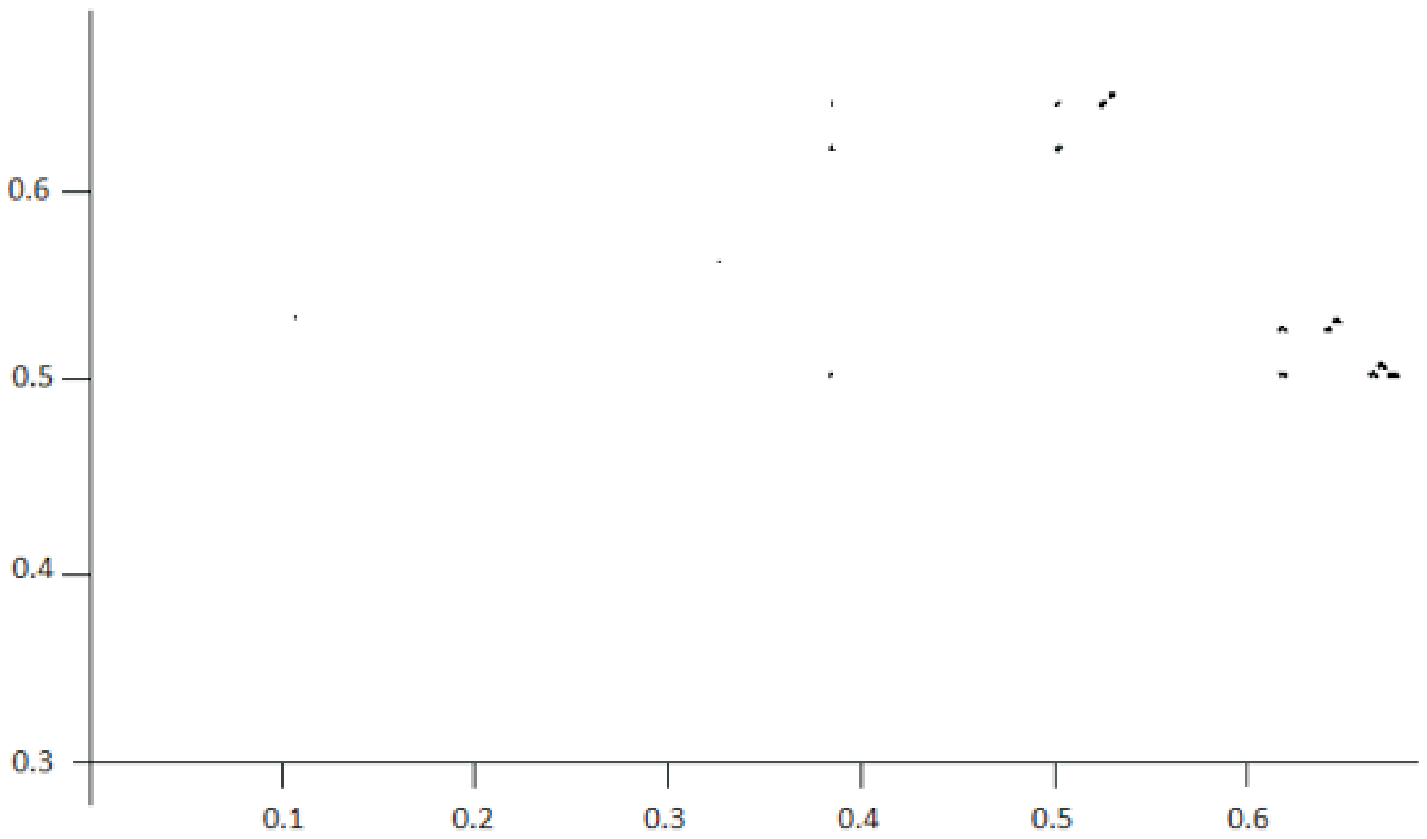}\\
Figure 5
\end{center}

\bigskip
\begin{center}
\includegraphics[width=13cm]{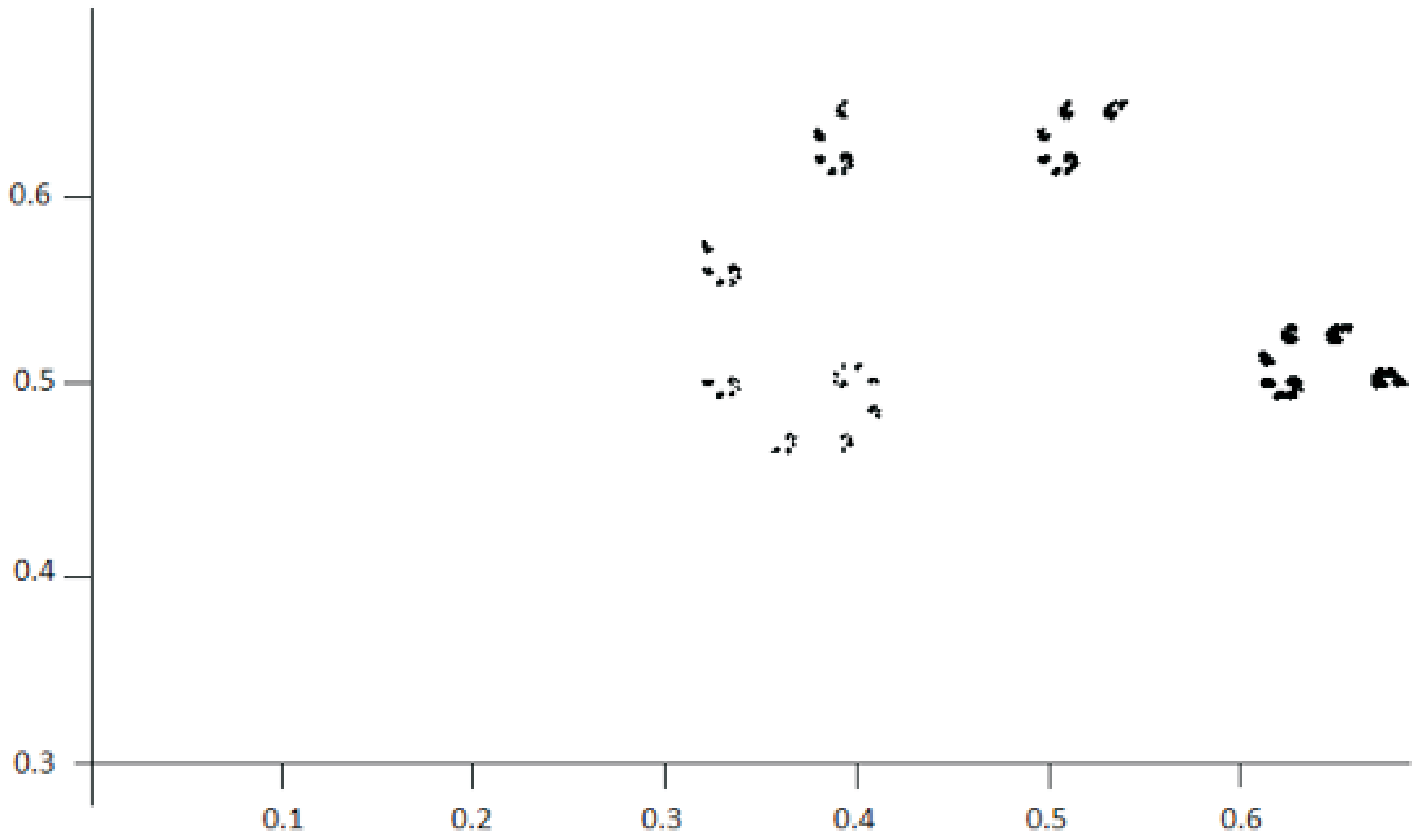}\\
Figure 6
\end{center}

\bigskip
 \textbf{Acknowledgements. }The authors are very grateful to the reviewers
whose extremely generous and valuable remarks and comments brought
substantial improvements to the paper.

\bigskip

\textbf{References}

\bigskip

[1] F. Browder, On the convergence of successive approximations for
nonlinear functional equations, Indag. Math., \textbf{30} (1968), 27--35.

[2] J.E. Hutchinson, Fractals and self similarity, Indiana Univ. Math. J., 
\textbf{30 }(1981), 713-747.

[3] P. Jaros, \L . Ma\'{s}lanka and F. Strobin, Algorithms generating images
of attractors of generalized iterated function systems, Numer. Algorithms, 
\textbf{73} (2016), 477-499.

[4] \L . Ma\'{s}lanka, F. Strobin, On generalized iterated function systems
defined on $l^{\infty }$-sum of a metric space, J. Math. Anal. Appl., 
\textbf{461} (2018), 1795-1832.

[5] J. Matkowski, Integrable solutions of functional equations, Diss. Math.,
127 (1975), 68 pp.

[6] R. Miculescu, Generalized iterated function systems with place dependent
probabilities, Acta Appl. Math., \textbf{130} (2014), 135-150.

[7] A. Mihail and R. Miculescu, Applications of Fixed Point Theorems in the
Theory of Generalized IFS, Fixed Point Theory Appl. Volume 2008, Article ID
312876, 11 pages doi: 10.1155/2008/312876.

[8] A. Mihail and R. Miculescu, A generalization of the Hutchinson measure,
Mediterr. J. Math., \textbf{6} (2009), 203--213.

[9] A. Mihail and R. Miculescu, Generalized IFSs on Noncompact Spaces, Fixed
Point Theory Appl. Volume 2010, Article ID 584215, 11 pages doi:
10.1155/2010/584215.

[10] E. Oliveira, The Ergodic Theorem for a new kind of attractor of a GIFS,
Chaos Solitons Fractals, \textbf{98} (2017), 63--71.

[11] E. Oliveira and F. Strobin, Fuzzy attractors appearing from GIFZS,
Fuzzy Sets Syst., \textbf{331} (2018), 131-156.

[12] N.A. Secelean, Invariant measure associated with a generalized
countable iterated function system, Mediterr. J. Math., \textbf{11} (2014),
361-372.

[13] N.A. Secelean, Generalized iterated function systems on the space $%
l^{\infty }(X)$, J. Math. Anal. Appl., \textbf{410} (2014), 847-858.

[14] F. Strobin, Attractors of generalized IFSs that are not attractors of
IFSs, J. Math. Anal. Appl., \textbf{422 }(2015), 99-108.

[15] F. Strobin and J. Swaczyna, On a certain generalisation of the iterated
function system, Bull. Aust. Math. Soc., \textbf{87} (2013), 37-54.

[16] F. Strobin and J. Swaczyna, A code space for a generalized IFS, Fixed
Point Theory, \textbf{17} (2016), 477-493.

\bigskip

{\small Radu MICULESCU}

{\small Faculty of Mathematics and Computer Science}

{\small Transilvania University of Bra\c{s}ov}

{\small Iuliu Maniu Street, nr. 50, 500091}

{\small Bra\c{s}ov, Romania}

{\small E-mail: radu.miculescu@unitbv.ro}

\bigskip

{\small Alexandru MIHAIL}

{\small Faculty of Mathematics and Computer Science}

{\small University of Bucharest, Romania}

{\small Academiei Street 14, 010014, Bucharest, Romania}

{\small E-mail: mihail\_alex@yahoo.com}

\bigskip

{\small Silviu-Aurelian URZICEANU}

{\small Faculty of Mathematics and Computer Science}

{\small University of Pite\c{s}ti, Romania}

{\small T\^{a}rgul din Vale 1, 110040, Pite\c{s}ti, Arge\c{s}, Romania}

{\small E-mail: fmi\_silviu@yahoo.com}

\end{document}